\documentclass[12pt]{amsart}
\newtheorem{theorem}{Theorem}[section]
\newtheorem{lemma}[theorem]{Lemma}
\newtheorem{proposition}[theorem]{Proposition}
\newtheorem{corollary}[theorem]{Corollary}
\theoremstyle{definition}
\newtheorem{definition}[theorem]{Definition}
\newtheorem{example}[theorem]{Example}
\theoremstyle{remark}
\newtheorem{remark}[theorem]{Remark}
\numberwithin{equation}{section}
          
\newcommand{\sdir}[1]{\rtimes_{#1}}
\def\R{{\mathbb R}}
\def\Z{{\mathbb Z}}

\def\L{{\mathfrak L}}

\def\SS{{\mathbb S}}
\def\T{{\mathbb T}}
\def\cE{{\mathcal E}} 
\def\cF{{\mathcal F}}
\def\til{\widetilde}
\def\defi{:=}

\def\id{\mathop{\rm id}\nolimits}
\def\Aut{\mathop{\rm Aut}\nolimits}
\def\SAut{\mathop{\rm SAut}\nolimits}

\def\hat{\widehat} \def\til{\widetilde}

\def\.{{\cdot}}
\def\&lt;{\langle}  \def\<{\langle}
\def\&gt;{\rangle}  \def\>{\rangle}
\def\ssk{\smallskip}       
\def\msk{\medskip}
\def\bsk{\bigskip}

\long\def\alert#1{\parindent2em\smallskip\hbox to\hsize%
{\hskip\parindent\vrule%
\vbox{\advance\hsize-2\parindent\hrule\smallskip\parindent.4\parindent%
\narrower\noindent#1\smallskip\hrule}\vrule\hfill}\smallskip\parindent0pt}

\newcommand{\gp}[1]{\langle#1\rangle}

\def\gen #1{\overline{\langle #1\rangle}}

\def\N{{\mathbb N}}             
\def\Q{{\mathbb Q}}             
\def\Saut{\mathop\mathrm{SAut}} 

\newcommand{\ignore}[1]{} 
\newcounter{alert}0

\newcommand{\na}{near abel\-ian}

\newcommand{\im}{in\-duc\-tive\-ly mo\-no\-the\-tic}

\def\.{{\cdot}}

\newcommand{\tqh}{to\-po\-logi\-cal\-ly qua\-si\-ha\-mil\-ton\-ian}

\newcommand{\lc}{loc\-al\-ly com\-pact}
\def\lca{\lc\ ab\-el\-ian}

\def\cC{{\mathcal S}\hskip-.5pt{\mathcal U}\hskip-.9pt{\mathcal B}}

\newcommand{\tlf}{{compactly ruled}}
\newcommand{\ccv}{\tlf}

\def\lead{\leaders\hbox to 1.5ex{\hss${.}$\hss}\hfill}
\def\arr{\hbox to 60pt{\rightarrowfill}}
\def\larr{\hbox to 60pt{\leftarrowfill}}

\newcommand{\psyl}[1]{$#1$-Sylow\,}
\newcommand{\sigcpt}{sigma-compact}

\usepackage{enumerate}
\usepackage{amssymb}
\usepackage[all]{xypic}
\usepackage{tikz} \usetikzlibrary{intersections}
\usepackage{tkz-euclide} \usetkzobj{all} 
\usepackage{comment}
\usepackage{color}

\begin{document}

\title[Locally Compact Near-Abelian Groups]
{A Study in Locally Compact Groups\\
\centerline{\sc ---Chabauty Space, Sylow Theory,} 
\centerline{\sc the Schur-Zassenhaus  Formalism,}
\centerline{\sc the Prime Graph for Near Abelian Groups}
}

\author[W. Herfort]{Wolfgang Herfort}
\address{
W. Herfort: Institute for Analysis and Scientific Computation,  \newline
\hspace*{5mm}Technische Universit\"at Wien, 
A-1040 Vienna, Austria} 
\email{wolfgang.herfort@tuwien.ac.at}

\author[K. H. Hofmann]{Karl H. Hofmann}
\address{
K.H. Hofmann: Fachbereich Mathematik,
Technische Universit\"at Darmstadt,\newline 
\hspace*{5mm}64289 Darmstadt, Germany} 
\email{hofmann@mathematik.tu-darmstadt.de}

\author[F. G. Russo]{Francesco G. Russo}
\address{
F.G. Russo: Department of Mathematics and Applied Mathematics,\newline 
\hspace*{5mm}University of Cape Town,
Rondebosch 7701, 
Cape Town, South Africa}
\email{francescog.russo@yahoo.com}

\dedicatory{Dedicated to Professor Herbert Heyer on the occasion of his 
eightieth birthday}

\begin{abstract}
The class of {\em \lc\ \na\ groups} is introduced and investigated
as a class of metabelian groups formalizing and applying the 
concept of scalar multiplication. 
 The structure of
\lc\ \na\ groups and its close connections to prime number theory  
are discussed and elucidated by graph theoretical
tools. These investigations require a thorough reviewing and 
extension to present circumstances of 
various aspects of the  general theory of locally compact groups
 such as 
\hfill\break
--the Chabauty  space of closed subgroups with its 
natural compact Hausdorff  topology,
\hfill\break
--a very general Sylow subgroup theory for periodic groups
including their Hall systems,
\hfill\break
--the scalar automorphisms of locally compact abelian groups, 
\hfill\break
--the theory of products of closed subgroups and their relation
to semidirect products, and 
\hfill\break
--inductively monothetic groups are introduced and
classified.
\hfill\break
As  applications, firstly, a complete  classification is given
of {\em \lc\ \tqh\ groups}, which has been initiated by F.~K\"um\-mich, 
and, secondly,  Yu.~Mukhin's  
classification of {\em \lc\ topologically modular groups} is retrieved and 
 further illuminated. 
\end{abstract}



\newtheorem{notation}[theorem]{Notation}
\newtheorem{remarks}[theorem]{Remarks}
\newtheorem{conjecture}[theorem]{Conjecture}

\maketitle

\section{Background}

In this survey we describe what may be called a structure theory 
of locally compact near-abelian groups. 
This attempt 
becomes clearer after we describe the historical development 
of researching  locally compact groups in a
broad sense.

Keys to our understanding will be the concept of 
an  inductively monothetic locally compact group, 
 the Chabauty space associated canonically 
with a locally compact group, and  a Sylow theory
for closed subgroups of periodic  locally compact groups reflecting typical
features of the Sylow theory of finite groups. The emergence of the role of
the prime numbers attached to the building blocks of the groups we consider
points into the direction of graph theory which we shall employ rather
extensively.   However,  as interesting as the structure
theory of near abelian groups itself is, we emphasise here,
that the machinery for developing
and detailing it should share center stage.

 \subsection{Some History} 
\qquad The structure theory of locally compact groups has a long
history.  One of its roots is {\sc David Hilbert}'s question of
1900 whether a locally euclidean topological group might
possibly support the introduction of a differentiable 
parametrisation such that the group operations are 
in fact differentiable.

A first affirmative answer was given for  {\em compact} locally
euclidean groups when they were found to be real matrix groups as 
a consequence of the 
foundational  work by {\sc Hermann Weyl} and his student 
{\sc Fritz Peter}
in 1927 on the representation theory of compact groups. 
A second step was achieved when the answer
was found to be ``yes'' for {\em commutative}  locally compact 
groups. This emerged out of the fundamental duality theorems
by {\sc Lev Semyonovich Pontryagin} (1934)
and {\sc Egbert van Kampen} (1937). This duality  forever determined the
structure and harmonic analysis of locally compact abelian
groups after it was widely read in the book of 1938 and 1953 by
{\sc Andr\'e Weil} on the integration in locally compact groups \cite{int}. 

\msk

A final positive answer to Hilbert's Fifth Problem
had to wait almost another
two decades when, 17 years after
 the Second World War,  the contributions of
{\sc Andrew Mattei Gleason, Deane Montgomery}, and {\sc Leo Zippin}
around 1952 provided the final affirmative answer
to {\sc Hilbert}'s problem. It led  almost at once to  the
fundamental insights of {\sc Hidehiko Yamabe} 1953, 
completing the pioneering work of {\sc Kenkichi Iwasawa} (1949)
providing  the fundamental structure of all those locally compact
groups $G$ which had a compact space $G/G_0$ of connected
components: Such groups were recognised as being approximated 
by quotient groups $G/N$ modulo arbitrarily small compact normal 
subgroups $N$ in such a fashion that each $G/N$ is a 
{\em Lie group}, that is, one of those groups  on which Hilbert 
had focussed in the fifth of his 23 influential  problems in 1900
and which {\sc Sophus Marius Lie} (1842-1899)
had invented together with an ingenious algebraisation 
method, long known nowadays under the name of Lie algebra
theory. (S. \cite{montzip}, \cite{Iwa3}) A special case
arises when all $G/N$ are discrete finite groups; in this
case $G$ is called {\em profinite}.

\msk

The solution of Hilbert's Fifth Problem
in the middle of last century opened up the access to
the structure theory of locally compact groups to the extent
they could be approximated by Lie groups, due to the rich Lie 
theory meanwhile developed in algebra, geometry, and functional
analysis. Recently, interest in 
{\sc Hilbert's Fifth} Problem was rekindled in he present
century under the influence of {\sc Terence Tao}
\cite{tao}. 

The quest for a solution to Hilbert's Fifth problem, at any rate,
led to one major direction in  the research of topological
groups: in focus was the class of groups $G$ 
approximable by Lie group quotients
$G/N$, and finally   $G$ itself needed no longer 
to  be  locally compact.
Such groups were called {\em pro-Lie groups} considered
for the sake of their own. (See \cite{hofmor-conn}, \cite{hofmor-almconn}.)
Their theory reached as far as {\em almost connected} locally compact
groups go (that is, those for which $G/G_0$ is compact),
including all compact
and all connected locally compact groups. But not further.
  
\ssk
Still, every locally compact group (and every pro-Lie group) $G$ has
canonically and functorially attached to it a (frequently infinite
dimensional) Lie algebra $\L(G)$ and therefore a cardinal $\ge0$ attached to
it, namely, its topological dimension $\mathrm{DIM}\, G$ 
(cf.\ e,g, \cite{hofmor}, 9.54, p.~498ff.). Indeed, we have the following
information right away:

\begin{proposition} \label{p:dim} Let $G$ be a topological group which is
\lc\ or a pro-Lie group. Then the following statements are equivalent:
\begin{enumerate}[\rm(1)]
\item Every connected component of $G$ is singleton, that is, $G$ is
      totally disconnected.

\item  $\mathrm{DIM}\, G=0$, that is, $G$ is zero-dimensional.
\end{enumerate}\end{proposition}

If $G$ is compact, then it is zero-dimensional if and only if it is 
profinite.
 
\bsk

For this survey it is important to see that, in the 20th century,
 there was a second trend on the 
study of locally compact
groups that is  equally significant even though it is opposite
to the concept of connectivity in topological groups. 
This trend is represented by the class of compact or locally
compact zero-dimensional groups.

Such groups were encountered in field theory 
at an early stage.  Indeed,  in  Galois theory  
the consideration of the appropriate infinite ascending 
family of finite Galois extensions and, finally, its  union 
would, dually, lead to an  inverse family of finite Galois groups 
and, in the end, to their projective limit. Thus was produced 
what  became known as a profinite group. The 
Galois group of the infinite field extension, 
equipped with the {\em Krull topology}, is thus a profinite group. 
One recognised soon
that profinite groups and
compact totally disconnected groups were one and the same 
mathematical object, expressed algebraically on the one hand and
topologically on the other. Comprehensive literature on this 
class of groups appeared 
much later than text books  on topological groups in which 
connected components
played a leading role. Just before the end of the 20th century
totally disconnected compact groups 
were the protagonists
of books simultaneously entitled ``Profinite Groups''
by {\sc John Stuart Wilson}  and  {\sc Luis Ribes} jointly with 
{\sc Pavel Zalesskii} in 1998, while {\sc George Willis} in
1994 laid the foundations of a general structure theory of
totally disconnected locally compact groups if
no additional algebraic information about them is available.

On the other hand, in the realm of locally compact abelian groups,
the completion of the field $\Q$ of rational numbers
with respect to any {\em nonarchimedian} valuation yields the
locally compact  $p$-adic fields $\Q_p$ as a totally disconnected
counterpart of the connected field $\R$ of real numbers. The fields $\Q_p$
and their  integral subrings
$\Z_p$ were basic building blocks of ever so many totally disconnected
groups, in  particular the linear groups over these field and
indeed all $p$-adic Lie groups which {\sc Nicolas Bourbaki}  judiciously
included in his comprehensive treatise on Lie groups. {\sc Bourbaki}'s
text on Lie groups formed the culmination and certainly the endpoint
of his encyclopaedic project extending over several decades.
 
The world of totally disconnected locally compact groups developed its
own existence, methods and philosophies, partly deriving from
finite group theory via approximation through the formation of
projective limits, partly through graph theory where the appropriate
automorphism theory provides the fitting representation theory, and 
partly also through the general impact of algebraic number theory
which in the text book literature is indicated by the books
of {\sc Helmut Hasse} since the thirties  and lastly, 1967,  by 
the ``Basic Number Theory'' of  {\sc Andr\'e Weil},
whose book on the integration on locally compact groups
of 1939/1953 had influenced the progress of harmonic analysis of
locally compact groups so much.

\bsk

We pointed out that each locally compact group $G$, irrespective of
any structural assumption has attached to it a (topological) Lie 
algebra  $\L(G)$ (and therefore a universal dimension).
The more recent interest in zero-dimensional locally compact groups $G$ 
has led to a new focus on another functorially attached  invariant,
namely, a compact Hausdorff space $\cC(G)$ consisting of all closed subgroups
of $G$ endowed with a suitable topology and now frequently called
the {\em Chabauty space} of $G$. 
This tool  is not exactly new, but has been widely utilised in 
applications  recently.  The Chabauty space is a special case of
what has been called the {\em hyperspace} of a compact (or locally compact)
space first introduced by {\sc Vietoris} (s.~\cite{engelking})
In topological algebra hyperspaces were used and described 
e.g.\ in \cite{bourbaki29}, \cite{hofmos}, \cite{waf}, and in
 all recent publications
where the name of Chabauty appears in the title (e.g.\ \cite{ham}, 
\cite{hamr}, \cite{hamhof}).
\msk

{\noindent\bf Notation.} We shall make use of notation coherent with the book of Hofmann and Morris,
\cite{hofmor}. The cyclic group of order $p^n$ is denoted by $\Z(p^n)$ and   $\Z(p^0)$ stands for the 
trivial abelian group $\{0\}$.

\section{Introductory Definitions and Results}

We now return to the second thrust of the study of locally compact
groups which is concerned with the research of $0$-dimensional 
groups.

\begin{definition} \label{d:periodic} \rm A topological group $G$
is called {\em periodic} if 
\begin{enumerate} [(i)]
\item $G$ is \lc\ and totally disconnected
\item $\gen g$ is compact for all $g\in G$.
\end{enumerate}
\end{definition} 

So a compact group is periodic if and only if it is profinite.
A very significant portion of the locally compact
groups considered here  will be {\em periodic} groups.
That is, we deal with totally disconnected \lc\  groups in which every 
element is contained in a profinite subgroup. 
We shall say
that a topological group $G$ is {\em compactly ruled} if it is the 
directed union of its compact open subgroups. 
{\em If $G$ is a a locally compact
solvable group in which every element is contained in a compact subgroup, then
it is compactly ruled.} 
The class of {\em compactly ruled groups} comprises both, the class of
{\em profinite groups} and the one of {\em locally finite} groups, i.e. groups
where every finite subset generates a finite subgroup only, see \cite{KW73}.
In many
cases we assume that the periodic groups we consider are
compactly ruled. These properties make them topologically special;
just how close they make our groups to profinite groups 
remains to be seen in the course of this survey.
A second significant property of the groups we study
is an algebraic one:
 they are solvable, indeed metabelian. Again, it is another
challenge  to discern
 just how close this
makes them to abelian groups.
The groups we study will be called {\em \na}.

\bsk

In order to offer a precise definition of the class of \lc\ groups
we  need one preliminary definition, extending a
very  familiar concept:

\begin{definition} \label{d:ind mon}\rm A topological group $G$ will
be called {\em monothetic} if $G=\gen g$ for some $g\in G$,
and {\em \im} if for every finite subset $F\subseteq G$ there is 
an element $g\in G$   such that 
$\gen F=\gen g$.
\end{definition}

We shall discuss and classify \im\ \lc\ 
groups in greater detail
later; but let us observe here right away a connected example
illustrating the two definitions:  Indeed let $\T=\R/\Z$ denote the
(additively written) circle group. Then 

{\em the 2-torus is monothetic but is not inductively monothetic, 
since $~(\frac 1 2\.\Z/\Z)^2\subseteq\T^2$ is finitely generated 
but is not monothetic.}

Yet in the domain of totally disconnected \lc\ groups 

{\em every $0$-dimensional monothetic
group is inductively monothetic.}

 In a periodic group, each monothetic
subgroup $\gen g$ is compact, equivalently, procyclic. 

Now we are prepared for a definition of 
the class of \lc\ groups whose details we shall consider here:

\begin{definition}\label{d:na defined}\rm
A topological  group $G$ is {\em \na}
provided it is locally compact and
contains a closed  abelian subgroup $A$ such that 
\begin{enumerate}
\item $G/A$ is an abelian \im\ group, and
\item every closed subgroup of $A$ is normal in $G$.
\end{enumerate} \end{definition}

The subgroup $A$  we shall call a 
{\em base} for $G$.

When we eventually collect   applications for this class of
locally compact groups, then we shall see that for instance
all locally compact groups in which two closed subgroups commute setwise
form a subclass of the class of \na\ groups and that the class of all 
locally compact groups in which the lattice of closed subgroups is modular
is likewise a subclass of the class of \na\ groups. 

\subsection{Some History of Near Abelian groups}

 In the world of discrete groups,
 \na\ groups historically appeared in a natural way
 when  {\sc K. Iwasawa}  attempted the classification
what is now known as {\em quasihamiltonian} and
{\em modular} groups as expounded in  the monograph by 
{\sc R. Schmidt}, \cite{schmidt}).
It was {\sc F. K\"ummich} (cf. \cite{kummich1}) who initiated 
in his dissertation  
written under the direction of {\sc Peter Plaumann} and in papers developed
from his thesis the study
of {\em topologically quasihamiltonian groups}. These are topological
groups such that $\overline{XY}=\overline{YX}$ is valid for any
closed subgroups $X$ and $Y$ of such a group. 
A bit later {\sc Yu. Mukhin} turned to investigating the 
class of locally compact 
topologically modular groups (cf. e.g. \cite{muk2,muk5}). 

\bsk

The properties that there be a closed normal abelian subgroup
$A$ of $G$ such that $G/A$ is inductively monothetic and such that
every closed subgroup of $A$ is normal in $G$ 
 suggest themselves  by the fact, proved by 
{\sc K. Iwasawa}, 
that discrete quasihamiltonian groups satisfy them. 
In a similar vein, {\sc Mukhin}, during his
work on classifying topologically modular groups, 
finds that these groups are all \na\ in our sense (see e.g. \cite{muk3}). 

An earlier article  by  {\sc K. H. Hofmann} 
and {\sc F. G. Russo},
was devoted to {\em classifying} compact $p$-groups that are topologically
quasihamiltonian (cf. \cite{hofrus}). 
The major result states that such groups are at the same
time topologically quasihamiltonian and \na\ with the
exception of $p=2$ in which case some sporadic \na\ groups
are topologically quasihamiltonian while the bulk of them are not.
This once again is evidence of the fact that is often quoted
by number theorists and group theorists alike that
$2$ is the oddest of all primes. 

\bsk

In linear algebra, a group $G$ of $(n+1)\times(n+1)$-matrices of
the form 
$$ \begin{pmatrix} r\.E_n&v\cr 0&1\cr\end{pmatrix},\quad  0<r\in\R,\ v\in\R^n$$
with the identity $E_n$ of GL$(n,\R)$
is a metabelian Lie group that has been called {\em almost abelian} 
(s.~\cite{cone}, p.~408, Example V.4.13). The subgroup $A$ of
all matrices with $r=1$ is isomorphic to $\R^n$ and every vector subspace
of $A$ is normal in $G$ and $G/A\cong\R$ is a one-dimensional Lie
group which is not \im, but we shall see that \im\ groups are
in some sense ``rank one'' group analogs.

In both cases we have a representation $\psi\colon G/A\to \Aut(A)$
such that $\psi(gA)(a)=gag^{-1}$ as an essential element of
structure.  In the near abelian case we shall say that $G$ 
is $A$-{\em nontrivial} if the image of $\psi$ has more than $2$ elements. 
 Whereas in the Lie group case, the structure of an almost abelian
Lie group $G$ is comparatively
simple, in the case of a group $G$ satisfying the conditions of 
Definition \ref{d:na defined} it is likely to be rather sophisticated as we
illustrate by a   result (s. \cite{BOOK}, Theorem 7.4) in which
$C_G(A)=\ker\psi$ denotes the centraliser $\{g\in G:(\forall a\in A)\, ag=ga\}$
of $A$ in $G$:

\begin{theorem} \label{th:na-nontrivial} 
{\rm(Structure Theorem I on Near Abelian Groups)} 
Let $G$ be an $A$-nontrivial  \na\ group.
Then 
\begin{enumerate}[{\rm (1)}]
\item $A$ is periodic.
\item $G$ is totally disconnected.
\item When $\psi(G/A)$ is compact or $A$ is an open subgroup,
      then $G$ has arbitrarily small compact open normal subgroups,
      that is, $G$ is {\em pro-discrete}.
\item $G$ itself is periodic if and only if $G/A$ is periodic if and only if
      $G/A$ is not isomorphic to a subgroup of the discrete group $\Q$ of rational
      numbers.
\item $C_G(A)$ is an abelian normal subgroup containing $A$ and is maximal
      for this property.
\end{enumerate}\end{theorem}

This shows that for our topic, {\em periodic} \lc\ groups play a significant 
role.

\section{Inductively Monothetic Groups}

A good understanding of near-abelian groups depends on 
a clear insight into the concept of  inductively monothetic groups. 
They were recently featured in \cite{hamhof}.

\msk

We must recall the  concept of a {\em local product} of a family
of topological groups which in the theory of locally compact groups mediates
between the idea of a Tychonoff product of compact groups and the idea
of a direct sum of a family of discrete groups; the principal applications
are in the domain of abelian groups, but the concept as such has nothing to
do with commutativity.

\begin{definition}\label{d:loc prod} \rm Let $(G_j)_{j\in J}$ be a family of 
locally compact groups and assume that for each $j\in J$ the group
$G_j$ contains a compact open subgroup $C_j$. Let $P$ be the   
subgroup of the cartesian product of the $G_j$ containing
exactly those $J$-tuples $(g_j)_{j\in J}$ of elements $g_j\in G_j$
for which the set $\{j\in J: g_j\notin C_j\}$ is finite. Then 
$P$ contains the cartesian product $C:=\prod_{j\in J} C_j$ which is
a compact topological group with respect to the Tychonoff topology.
The group $P$ has a unique group topology with respect to which $C$
is an open subgroup. Now 
the {\em local product} of the family $((G_j,C_j))_{j\in J}$ is the group $P$
with this topology, and it is denoted by
$$P=\prod_{j\in J}^{\rm loc}(G_j,C_j).$$
\end{definition}

Let us note that the local product is a \lc\ group with the
compact open subgroup $\prod_{i\in J}C_p$. While the full product
$\prod_{j\in J}G_j$ has its own product topology we note that
in general the local product topology on $P$ is properly
finer than the subgroup topology. The concept of the local product
was introduced and its duality theory in the commutative situation was   
studied by {\sc J. Braconnier} in \cite{braconnier}.
For us local products play a role most frequently with $J$ being
the set $\pi$ of all prime numbers. This is well illustrated 
by the following key result on periodic \lca\ groups where we note,
that for a \lca\ group $G$ and each prime $p$,
we have a unique characteristic subgroup $G_p$ containing all elements
$g$ for which $\gen g$ is a profinite $p$-group; $G_p$ is called
the $p$-{\em primary component} or the $p$-{\em Sylow subgroup} of $G$.

\begin{theorem}\label{th:braconnier}\quad {\rm (J. Braconnier)}
Let $G$ be a  
periodic \lc\ abel\-ian group and $C$ any compact open subgroup of $G$. Then
$G$ is isomorphic to the local product 
$$\prod_{p\in\pi}^{\rm loc} (G_p,C_p).\eqno{\rm(LP)}$$
\end{theorem}

The following remark is useful for us as a consequence of the fact that
any compact $p$-group $C_p$ is a $\Z_p$-module and any prime $q\ne p$ is a
unit in $\Z_p$, whence $C_p$ is divisible by $n\in \N$ with $(n,p)=1$:

\begin{remark} \label{r:divi} A periodic \lca\ group $G$ is divisible
iff all $p$-Sylow subgroups $G_p$ are divisible.
\end{remark}

\bsk 

The structure of a locally compact {\em monothetic} group $G$ is familiar
to workers in the area: It is either isomorphic to the discrete group $\Z$
of integers or is compact 
 (Weil's Lemma, s.\  
e.g.\ \cite{hofmor}, Proposition 7.43, p.348.).
A compact abelian group is known if its discrete Pontryagin dual is known.
A compact abelian  group $G$ is monothetic if and only if there is
a morphism $f\colon \Z\to G$ of locally compact groups with dense image,
that is, iff there is an injection of the discrete group $\hat G$ into
the character group $\T$, that is, $\hat G$ is isomorphic to
a subgroup of $\Q^{(2^{\aleph_0})}\oplus\bigoplus_{p\in\pi}\Z(p^\infty)$.
(Here $\Z(p^\infty)$, as usual, is the {\it Pr\"ufer group} 
$\bigcup_{n\in\N}\frac 1{p^n}\Z/\Z\subseteq \T$.) 
Whenever $G$ is zero-dimensional, things simplify dramatically:

\begin{proposition} \label{l:0-monoth} A compact zero-dimensional abelian
group $G$ is monothetic if it is isomorphic to 
$\prod_{p\in\pi}G_p$ where the $p$-factor $G_p$ is either $\Z(p^m)$
for some $m\in \N_0=\{0,1,2,\dots\}$ or $\Z_p$, the additive group
of the ring of $p$-adic numbers.
\end{proposition} 

Let us proceed to \im\ \lc\ groups. For periodic \im\ groups it
 is convenient to introduce some special terminology. From Braconnier's
Theorem \ref{th:braconnier} we know that every periodic \lca\ group $G$,
for any given compact open subgroup $C\subseteq G$,
is (isomorphic to) the local product 
$$\prod_{p\in \pi}(G_p,C_p)\eqno(*)$$ of
its $p$-Sylow subgroups.  

\begin{definition}\label{d:special defs} \rm
 A topological group $G$ is called
$\Pi$-{\em procyclic}, if it is a periodic \lca\ group and
each $p$-Sylow subgroup $G_p$ is either a finite cyclic
$p$-group (possibly singleton) or $\Z_p$, that is, $G_p$ is 
$p$-procyclic.
\end{definition}

Now we can formulate the classification of \im\ \lc\ groups.

\begin{theorem}\label{th:thoughtsE} 
{\rm(Classification Theorem of Inductively Monothetic 
Groups)}
 Let $G$ be an inductively monothetic \lc\ group. 
Then  $G$ is either
\begin{enumerate}[\rm(a)]
\item a $1$-dimensional compact connected abelian group, or
\item a subgroup of the discrete group $\Q$, or
\item a periodic \lca\ group such that $G_p$ is isomorphic to 
$\Q_p$, or  
$\Z(p^\infty)$, or
$\Z_p$, or
$\Z(p^{n_p})$ for some $n_p\in\N_0=\{0,1,2,\dots\}$.
\end{enumerate}
\end{theorem}

All inductively monothetic groups are sigma-compact, i.e.\ are countable
unions of compact subsets.

The groups of connected type in
Theorem \ref{th:thoughtsE} 
are monothetic; other types may or may not 
be monothetic. The periodic \im\ groups $G$ of part (b) require special attention.
First we divide the set $\pi$ of all prime numbers into disjoint sets
\begin{itemize}
\item $\pi_A=\{p\in\pi: G_p\cong \Q_p\}$,
\item $\pi_B=\{p\in\pi: G_p\cong \Z(p^\infty)\}$,
\item $\pi_C=\{p\in\pi: G_p\cong \Z_p\}$,
\item $\pi_D=\{p\in\pi: (\exists n\in\N_0)\ G_p\cong\Z(p^n)\}$.
\end{itemize}
Now we fix a compact open subgroup $C$ of $G$ and
identify $G$ with the local product 
$\prod_{p\in \pi}(G_p,C_p)$, further we define two closed 
characteristic subgroups as 
$$\begin{matrix}  
 D&\defi \prod_{p\in \pi_A\cup \pi_B}^{\rm loc} (G_p,C_p),\\
 P&\defi \prod_{p\in \pi_C\cup \pi_D}^{\rm loc} (G_p,C_p),\\
\end{matrix}$$
and notice that $G=D\oplus P$. Both subgroups $D$ and $G$
are characteristic, and we notice that in view of 
Remark \ref{r:divi} $D$ is the unique
largest divisible subgroup of $G$

\begin{theorem}\label{th:thoughtsE2} 
{\rm(Classification  of Inductively Monothetic
Groups, continued)}
 Let $G$ be a periodic inductively monothetic \lc\ group. 
Then  $G$ is the direct topological and algebraic sum $D\oplus P$
of two characteristic closed subgroups of which $D$ is 
the largest divisible subgroup of $G$ and $P$ is the unique largest
$\Pi$-procyclic subgroup according to  {\rm Definition \ref{d:special defs}}.  
\end{theorem}

\bsk

We apply this information to the structure theory of a
near abelian periodic group $G$ with base $A$. Then 
 $G/A=D\oplus P$ as in 
Theorem \ref{th:thoughtsE2}; let $G_D$, respectively $G_P$ 
denote the full inverse images for the quotient morphism $G\to G/A$.
Now $G_D$ and $G_P$ are closed normal subgroups
such that $G=G_DG_P$ and $G_D\cap G_P=A$, and 
we have  $G_D\subseteq C_G(A)$ (see \cite{BOOK}, Theorem 7.6).

\begin{theorem} \label{th:na-nontrivial2} 
{\rm(Structure Theorem II on Near Abelian Groups)} 
Let $G$ be a periodic near abelian locally compact group with a base $A$
such that $G$ is 
  $A$-nontrivial. Then $A\subseteq G_D\subseteq C_G(A)$ 
 where $G_D$ is a normal abelian subgroup such that
 $G/G_D\cong G_P/A$ is $\Pi$-procyclic.
\end{theorem}

This portion of the basic structure theory of near abelian groups
in the periodic situation will allow us to concentrate largely on the case
that the factor group $G/A$ is $\Pi$-procyclic.

\section{Factorisation and Scaling}

We begin with a definition elaborating 
the definition of near abelian groups.

\begin{definition} \rm Let $G$ be a near abelian locally compact group
with a base $A$. A closed subgroup $H$ is called a {\em scaling
subgroup} for $A$ if 
\begin{enumerate}[\rm(i)]
\item $H$ is inductively monothetic, and
\item $G=AH$.
\end{enumerate}\end{definition}

\begin{example} \label{ex: maier} \rm

\quad There exists a (discrete) abelian group $G$ with a subgroup $A$
which is not a direct summand and which has the following properties:
$A$ is the torsion subgroup of $G$ of the form 
$A\cong \bigoplus_{n\in\N} \Z(p)$ and 
$G/A\cong \bigcup_{n\in\N}\frac 1{2\.3\cdots p_n}\Z\subseteq\Q$. 
\end{example}  

\begin{example} \label{ex: russo}
 $G\cong \bigoplus_{n\in\N}\Z(p^n)$ and there is a subgroup $A$ which 
is not a direct summand such that
$G/A\cong \Z(p^\infty)$.  
\end{example}

In Example \ref{ex: maier}, the group $G$ is a subgroup of $\R/\Z$ and is 
a  construction due to D. Maier \cite{maier}. Example \ref{ex: russo}
is inspired by Example $\nabla$ in 
Theorem A1.32, p.~686 of \cite{hofmor}. 

These examples show that there are  obstructions to a
very {\em general} result
asserting the existence of a scaling group for near abelian groups
$G$ with bases $A$.

\bsk

A scaling group $H$, whenever it exists, is a {\em supplement} for  
$A$ in $G$ but not in general a semidirect {\em complement}. How
far a supplement is from being a complement can be clarified under fairly
general circumstances;  we illustrate that in the following proposition.

\begin{proposition}\label{p: mayer-vietoris} Let $G$ be a 
\lc\ group with a closed normal subgroup $A$ and a closed 
\sigcpt\ subgroup $H$ containing a compact open subgroup and
 satisfying $G=AH$. The inner automorphisms define a morphism
$\alpha\colon H\to\Aut(A)$ by $\alpha(h)(a)=hah^{-1}$. 
Then we have the following conclusions:
\begin{enumerate}[\rm(i)]
\item The semidirect
product $A\sdir{\alpha}H$ is a  \lc\ group and the function 
$\mu\colon A\sdir{\alpha} H\to G$, $\mu(a,h)=ah$,
   is a quotient  morphism with kernel 
$\{(h^{-1},h):h\in A\cap H\}$ isomorphic to $A\cap H$, 
mapping both $A$ and $H$ faithfully.

\item The factor group $G/(A\cap H)$ is a semidirect product 
of $A/(A\cap H)$ and $H/(A\cap H)$ and
the composition 
$$A\rtimes_{\alpha}H\to G\to G/(A\cap H)$$
is equivalent to the natural quotient morphism 
$$A\rtimes_{\alpha}H\to\frac{A}{A\cap H}\rtimes\frac{H}{A\cap H}$$
with kernel $(A\cap H)\times (A\cap H)$.
\end{enumerate}
\end{proposition} 

Notice that a scaling subgroup $H$ of
a near abelian group is sigma-compact and has a compact open
subgroup, so that the proposition applies in its entirety to near abelian
locally compact groups. The typical {\em ``sandwich situation''
$$ A\rtimes H \to AH\to\frac{A}{A\cap H}\rtimes\frac{H}{A\cap H}$$ 
is also observed in significant ways in the structure theory
of compact groups (see \cite{hofmor}, e.g.~Corollary 6.75 ff.).} 

\bsk
So one of the most pressing questions of the structure 
theory of near abelian \lc\ groups is the following:

\msk
\noindent{\bf Problem 1.} {\em Under which conditions does a
locally compact group $G$ with a normal subgroup $A$ such that
$G/A$ is inductively monothetic 
contain a closed inductively monothetic subgroup $H$ such that $G=AH$?}
\bsk

If $G/A$ is in fact monothetic, then
 the answer is affirmative and
easy. In more general circumstances
we have the following theorems giving a partial answer to
Problem 1. 

\begin{theorem} \label{th:scaling1} 
Let $G$ be a  \lc\ group 
with a {\em compact} normal subgroup $A$ such that $G/A$ is $\Pi$-procyclic.
Then $G$ contain a $\Pi$-procyclic subgroup $H_\Pi$ such that $G=AH_\Pi$.
\end{theorem}

\begin{theorem} \label{th:scaling2} Let $G$ be a locally compact
group with a {\em compact open} normal subgroup $A$ such that 
$G/A$ is isomorphic to an infinite subgroup of the  group $\Q$.
Then $G$ contains a discrete subgroup $H\cong G/A$
such that $G$ is a semidirect product $AH\cong A\sdir{}H$.   
\end{theorem}

It would be highly desirable to have such  theorems without the hypothesis that
$A$ be {\em compact}. The proofs of these theorems 
(see \cite{BOOK}, Theorem 5.23ff.)  make essential use of the compact 
Hausdorff Vietoris-Chabauty space
$\cC(G)$ which is attached to {\em every} locally compact group
as a  general invariant. 

\ssk

As long as this approach requires the compactness
of $A$ the following theorem may be considered
as fundamental for the structure theory of near abelian 
locally compact groups:

\begin{theorem} \label{th:scaling3} Let $G$ be a locally compact near abelian
 group with a base $A$ such that $G$ is $A$-nontrivial and 
$G/A$ is $\Pi$-procyclic. Then $G$ contains a $\Pi$-procyclic
scaling subgroup $H_\Pi$ for $A$ with $G=AH_{\Pi}$.
\end{theorem}

For a proof see \cite{BOOK}, Theorem 11.3. 
The proof requires a wide spectrum of parts of our general
structure theory of near abelian groups.
In particular, at the root of this existence theorem is the
Chabauty space $\cC(G)$ of the group $G$ which we mentioned
earlier. 
This theorem and Proposition  \ref{th:na-nontrivial2} now yield the
following theorem (see \cite{BOOK}, Theorem 11.6). 

\begin{corollary} \label{c:ssg} For every periodic \lc\ \na\ group $G$ with
base $A$ there exists a $\Pi$-procyclic closed subgroup $H_\Pi$
such that $G=G_DH_\Pi$ for the abelian normal subgroup $G_D$ 
with $A\subseteq G_D\subseteq C_G(A)$.
\end{corollary}

(See \cite{BOOK}, Theorem 6.3(iv,v) and their proofs.)

In particular, Proposition \ref{p: mayer-vietoris} then shows us that
we have
\begin{corollary} \label{c:ssg2} Every periodic \lc\ \na\ group 
$G$ is a quotient of $G_D\sdir{}H_\Pi$ modulo a subgroup isomorphic
to $A\cap H_\Pi$.
\end{corollary}

Recall that $Z(G)$ denotes the center of a group $G$.

\begin{corollary} \label{c:ssg3} For every periodic \lc\ \na\ group $G$ with
base $A$ we have $C_G(A)=AZ(G)$ and $C_G(A)\cap H_\Pi\subseteq Z(G)$, that is
$AZ(G) \cap H_\Pi= Z(G)\cap H_\Pi$.
\end{corollary}  

(See \cite{BOOK}, Theorem 6.3(iv,v) and its proof.)

\bsk

The following theorem then is rather definitive on the factorisation of
a periodic near abelian locally compact group and may be considered
as one of the main theorems on their structure.

\begin{theorem} \label{c:ssg4} {\rm(Structure Theorem III on Periodic
Near Abel\-ian Groups)} 
For every periodic \lc\ \na\ group $G$ with
a base $A$ such that $G$ is not $A$-trivial, we have
$$G=AZ(G)H$$ with a $\Pi$-procyclic scaling group $H$ for $A$ in $G_P$,
where $AZ(G)\cap H=Z(G)\cap H$. 
\end{theorem}  
 
For information as to which closed subgroups $A^*\subseteq AZ(G)$
containing $A$ may still be taken as base subgroups, see 
\cite{BOOK}. Theorem 10.32. The role of the center $Z(G)$ 
in $C_G(A)=AZ(G)$---a locally compact abelian group we know to
be a local product of its $p$-primary components $A_pZ(G)_p$---
is still a bit mysterious; more information will be forthcoming
in Theorem \ref{th:Apq} below.

\section{The Sylow Theory of  Periodic Groups}

Sylow theory, i.e., existence and conjugacy of maximal $p$-subgroups,
and, more generally, 
of maximal $\sigma$-subgroups where $\sigma$ is a set of primes,
is available for profinite groups (see \cite{wilson,ribes-zalesskii}). Several attempts 
have been made in order to generalise Sylow theory to non-compact and
locally compact groups, see e.g. the survey from 1964 by \v Carin, \cite{Car64}, or,
more recently, Platonov in \cite{Pla66} and Reid in \cite{Reid13}.
Here we shall focus on our class of 
{\em compactly ruled groups}. If the topology on  a compactly ruled
 group is discrete,
the group is {\em locally finite}, i.e., every finite subset generates a finite
subgroup. Then our Sylow theory reduces to the one presented in the
book of Kegel and Wehrfritz, see \cite{KW73}. 

In each locally compact {\em periodic} group $G$ the
concept of a $p$-group can be defined meaningfully.
Indeed if $g\in G$, then $M\defi\gen g$ is a zero-dimensional monothetic
compact group, and thus
$$ M\cong\prod_{p\in \pi}M_p,$$
where for each prime $p$ the $p$-primary component $M_p$ 
is either $\cong \Z_p$ or $\Z(p^{n_p})$ for some $n_p=0,1,2,\dots$.

It is practical to generalise the concept of a $p$-element:
For each subset $\sigma\subseteq \pi$, an  element $g\in G$ is
 called a $\sigma$-element if $\gen g=\prod_{p\in\sigma}M_p$; if $\sigma
=\{p\}$, then $g$ is called a $p$-element.
The group $G$ is a $\sigma$-group, if all of its elements are $\sigma$-elements.
A subgroup $S$ is called a $\sigma$-Sylow subgroup of $G$ if it is 
a maximal element in the set of $\sigma$-subgroups.
A simple application of  Zorn's Lemma shows that every $\sigma$-element
is contained in a $\sigma$-Sylow subgroup. We record:
\begin{lemma} \label{l:closure} {\rm(The Closure Lemma)} Let $G$ be
any locally compact totally disconnected group. Then for any subset
$\sigma\subseteq\pi$, then set $G_\sigma$ of all $\sigma$-elements of
$G$ is closed in $G$.
\end{lemma}

Let us look to some traditional splitting theorems that still work in
the general background of periodic locally compact groups.

\subsection{The Schur-Zassenhaus Splitting} 
The splitting of finite groups into pro\-ducts
of subgroups of relatively prime orders can be generalized to the locally compact
setting up to a point, as we show in the following. For locally finite groups
the results to be discussed are well known, see e.g. \cite{KW73}. 
They also relate to work of the
second author, see \cite{H64,hofmos}.

\begin{proposition} \label{p:schz} Let $N$ be a closed
subgroup of a locally compact periodic group $G$ and assume
$N\subseteq G_\sigma$. Then the following conditions are equivalent:
\begin{enumerate}[$(1)$]
\item $N$ is a normal Sylow subgroup.
\item $N=G_\sigma$.
\item $N$ is normal and $G/N$ contains no $p$-element with $p\in\sigma$. 
\end{enumerate}
\end{proposition} 

\begin{definition}\label{d:schzc}\rm
Let $G$ be a \lc\ periodic group and $N$ a closed  subgroup.
We say that $N$ satisfies the {\em Schur-Zassenhaus Condition} 
if and only if it satisfies the equivalent conditions of 
Proposition \ref{p:schz} for $\sigma=\pi(N)$. 
\end{definition}

\begin{theorem}\label{th:schur-zassenhaus} \textrm{(Schur-Zassenhaus 
Theorem)}
Let $G$ be a periodic group and $N$ a closed  subgroup
satisfying the following two conditions:
\begin{enumerate}[\rm(1)]
\item $N$ satisfies the Schur-Zassenhaus Condition.
\item  $G/N$ is a directed countable union of compact subgroups.
\end{enumerate}
Then the following
conclusions hold:
\begin{enumerate}[\rm(i)]
\item  $N$ possesses a complement $H$ in $G$.
\item  Let $K$ be a closed subgroup of $G$ such that $K\cap N=\{1\}$
and assume that $G/N$ is compact.
Then there is a $g\in G$ such that $gKg^{-1}=H$.
\end{enumerate}
\end{theorem}

It should be remarked, that for solvable  groups (such as 
\na\ groups) the periodic groups are always directed unions
of their open compact subgroups. In such a situation condition
(2) simply means that $G/N$ is sigma-compact.

The Schur-Zassenhaus configuration in the locally compact environment
is  delicate, since, in general problems arise with the product of
a closed normal subgroup and a closed subgroup; such a product need not be 
closed, in general.  

Still we do have theorems like the following:

\begin{theorem}\label{t:normal-pi-Sylsgp}
Let $N$ be a normal $\sigma$-Sylow subgroup of a \lc\ periodic group $G$. Then
a $(\pi\setminus\sigma)$-Sylow subgroup $H$ of $G$ exists such that $NH$ is an open
and hence closed subgroup.  
Moreover, if $H$ is \emph{any} $(\pi\setminus\sigma)$-Sylow subgroup of $G$, then
$NH$ is closed in $G$ and $H$ is a complement 
of $N$ in $NH$, that is, $NH=N\rtimes H$.
\end{theorem}

\bsk
\subsection{Sylow Subgroups Commuting Pairwise}
Let $p\in \pi$ denote any prime and $p'\defi \pi\setminus\{p\}$. 

Then we have the following result:
\begin{lemma}  For a compactly ruled group $G$ and a prime number $p$,
the following conditions are
equivalent:
\begin{enumerate}[\rm(1)]
\item $[G_p,G_{p'}]=\{1\}$.
\item Both $G_p$ and $G_{p'}$ are subgroups, and $G=G_p\times G_{p'}$.
\item There is a unique projection ${\rm pr}_p\colon G\to G_p$ with
kernel $G_{p'}$.
\end{enumerate}\end{lemma}

\begin{definition} \rm For a periodic locally compact group $G$ we write
$$\nu(G)=\{p\in\pi: [G_p,G_{p'}]=\{1\}\}.$$
\end{definition}

We have found the following structure theorem very useful in the context
of near abelian groups generalising the well-known fact that a
{\em pronilpotent group} is the cartesian product of its Sylow subgroups 
(cf.~\cite{ribes-zalesskii}):

\begin{theorem} \label{th:nu-theorem} In a compactly ruled \lc\ group
$G$, the set $G_{\nu(G)'}$ of $\alpha$-elements with 
$\alpha\cap \nu(G)=\emptyset$ 
is a closed normal subgroup, and all $p$-Sylow subgroups for 
$p\in \nu(G)$  are normal subgroups. Moreover,
$$G\cong  G_{\nu(G)'}\times\prod_{p\in\nu(G)}^{\rm loc}(G_p,U_p)$$
for a suitable family of compact open subgroups $U_p\subseteq G_p$
as $p$ ranges through $\nu(G)$. 
\end{theorem}
 
\bsk
\subsection{The Internal Structure of Sylow Subgroups of Near Abelian \break Groups}

For periodic \na\ groups, to which we can apply a Sylow theory meaningfully
we  assume that $G$ is a periodic \na\ \lc\ group such that 
$G$ is nontrivial for a base $A$.
 
\begin{theorem} \label{th:period-sylow} 
Let $G$ be a periodic \na\ group and $A$ a base for which $G$ is 
$A$-nontrivial and which satisfies $A=C_G(A)$. 
Then, for every set $\sigma$ of prime numbers there 
is a \psyl\sigma subgroup $S_\sigma$.
Fix a \psyl\sigma subgroup $S_\sigma$ of $G$. Then
\begin{enumerate}[\rm(i)]
\item $S_\sigma\cap A$ is the $\sigma$-primary component $A_\sigma$ of $A$, equivalently,
the $\sigma$-Sylow subgroup  of $A$. Moreover, 
\item $S_\sigma/A_\sigma\cong S_\sigma A/A =(G/A)_\sigma$.
\item If $(G/C_G(A))_\sigma\cong H/(H\cap C_G(A))$ is compact, then
      any two \psyl\sigma\ subgroups of $G$ are conjugate.
\item $S_\sigma=C_G(A)_\sigma H_\sigma=A_\sigma Z(G)_\sigma H_\sigma$, where $H$ is  as in
{\rm Corollary \ref{c:ssg}}. 
\end{enumerate}
\end{theorem}
 
(See \cite{BOOK}, Theorem 10.1.)
The case that $\sigma=\{p\}$ is an important special case.
Let us note that it may happen, $S_p\subseteq A$, in which case we
have $H_p=\{1\}$.

\section{Scalar Automorphisms}

Among the methods we are using, the specification of
scalar automorphisms of a periodic \lca\ group is prominent.
Every locally compact abelian $p$-group $A$ is a natural $\Z_p$-module,
and in the case of a {\em periodic} \lc\ abelian group
$$A = \prod_{p\in\pi}^{\rm loc}(A_p,C_p) \eqno{\rm(LQ)} $$
it is a natural
$$ \til \Z= \prod_{p\in\pi} \Z_p$$
module by componentwise scalar multiplication 
$$z\.g=(z_p)_{p\in\pi}\.(g_p)_{p\in\pi}=(z_p\.g_p)_{p\in\pi}$$
The compact ring $\til\Z$ is the profinite compactification of the
ring $\Z$ of integers. 

\begin{lemma} \label{l:scalar} {\rm(The Scalar Morphism Lemma)} 
For a continuous automorphism $\alpha$ of a periodic  \lca\ group $G$
the following conditions are equivalent:
\begin{enumerate}[{\rm(1)}]
\item $\alpha(H)\subseteq H$ for all closed subgroups $H$ of $G$.
\item $\alpha(\gen g)\subseteq \gen g$ for all $g\in G$.
\item $\alpha(g)\in\gen g$ for all $g\in G$.
\item There is an $r\in\til\Z$ 
      such that $\alpha(a)=r\.a$ for all $g\in G$.
\end{enumerate}
\end{lemma}

We note in passing that the first three conditions are equivalent
in any locally compact group.

\begin{definition} \rm An automorphism $\alpha\in \Aut(A)$ 
of a periodic lo\-cal\-ly com\-pact abel\-ian group $A$ is called a {\em scalar automorphism}.
The group of all scalar automorphisms is written $\SAut(A)$.
\end{definition}

For $r\in\til\Z^\times$, the group of invertible elements of $\til\Z$,
 we denote the function $a\mapsto r\.a:A\to A$
by  $\mu_r\in\SAut(A)$.
 
\begin{proposition} \label{p:scal-auto} Let $A$ be a periodic \lca\ 
group. Then

\begin{enumerate}[\rm(i)]
\item $r\mapsto \mu_r: \til\Z^\times\to\SAut(A)$
is a quotient morphism of compact groups. 
In particular, $\SAut(A)$ is a profinite group and thus does not
contain any nondegenerate divisible subgroups.
\item The following conditions are equivalent:
\begin{enumerate}[\rm(a)]
\item $\SAut(A)=\{\id_A,-\id_A\}$,
\item The exponent of $A$ is $2$, $3$, or $4$.
\end{enumerate}
In particular, $A$ has exponent $2$ if and only if $-\id_A=\id_A$.
\end{enumerate}
\end{proposition}

In the process of these discussions, we recover in our framework the
following theorem of {\sc Mukhin} \cite{muk3}:

\begin{theorem} \label{th:muk3-th2}
Let $A$ be a \lca\  group written additively.
\begin{enumerate}[\rm(a)]
\item If $A$ is not periodic,  then $\Saut(A)=\{\id, -\id\}$. 
\item If $A$ is periodic, then $\Saut(A){=}\prod_p\Saut(A_p)$, 
              where $\Saut(A_p)$ may be identified with the group
               of units of  the ring of scalars of $A_p$:
              $$\begin{cases} \Z_p, &\mbox{if the exponent of $A_p$ is infinite}, \\
              \Z_p/p^m\Z_p\cong \Z(p^m), &\mbox{for suitable $m$ otherwise}.
              \end{cases}$$       
\item In particular, $\Saut(A)$ is a 
              homomorphic image of $\til\Z^\times$.
\item An automorphism $\alpha$ is in $\SAut(A)$ iff there is
              a unit $z\in\til\Z^\times$ such that
  
              \centerline{%
              \quad $(\forall g\in G)\, \alpha(g)=z\.g=\prod_p z_p\.g_p$ for 
              $z=\prod_pz_p$, $g=\prod_pg_p$.}
\end{enumerate}
\end{theorem}

The significance of Mukhin's Theorem for the structure theory of 
\na\ groups is visible in the very Definition \ref{d:na defined}
via Theorem \ref{c:ssg4}.
Indeed if $G$ is a \na\ \lc\ group with a base $A$, then the
inner automorphisms of $G$ induce a faithful action of the 
$\Pi$-procyclic  factor group $G/C_G(A)\cong H/(H\cap Z(G))$  upon the
base $A$.  
So  $\SAut(A)$ is a quotient of $H$
and therefore is $\Pi$-procyclic.    

The structure of $G$ is largely determined by the structure
of $\SAut(A)$ and therefore by the group ${\til\Z}^\times$ of units 
of $\til \Z$.

\subsection{The Group of Units of the Profinite Compactification of 
the Ring of Integers and its Prime Graph}

The group $\til\Z^\times$ is more complex than it appears at first.
Its Sylow theory or primary decomposition is best understood in
graph theoretical terms. The same graph theory turns out to be almost
indispensable for dealing with the Sylow structure of near abelian
groups in general. The graphs that we use are all subgraphs of 
a ``universal'' graph (which we also call the ``master graph'')
and which is used precisely to describe the Sylow theory of 
$\til\Z^\times$. We discuss it in the following.

\bsk

A {\em bipartite graph} consists of two disjoint sets $U$ and $V$ and 
a binary relation $E\subseteq (U\cup V)^2$ such that $(u,v)\in E$ 
implies $u\in U$ and \index{graph!bipartite}
$v\in V$. The elements of $U\cup V$ are called {\em vertices}
and the elements of $E$ are called {\em edges}. Any  triple
$(U,V,E)$ of this type is called a {\em bipartite graph}.  

\msk

In the following we construct a special bipartite graph 

\msk

\def\Gr{\mathcal G}

$\Gr=(U,V,E)$ with $U, V\subseteq \N\times \{0,1\}$ as follows:

\begin{definition} \label{d:bipartite} \rm
Define
$U=\N\times\{1\}$, $V=\N\times\{0\}$.

\ssk

Let $n\mapsto p_n$ be the unique order preserving bijection of
$\N$ onto the set $\pi$ of prime numbers.
On $\pi$ we consider the binary relation
$$ T=\{(p,q)\in\pi\times\pi: q=p\mbox{ or } p|(q-1)\}.\leqno(1)$$
Let $E\subseteq (\N\times\{0,1\})^2$ be defined as follows
$$E=\{((m,1),(n,0)): (p_m,q_n)\in T\}.\leqno(2)$$
If $e=((m,1),(n,0))\in E$ is an edge we shall use the following
notation for the prime numbers associated with $e$:
$$ p_e=p_m\qquad q_e=q_n.$$

\end{definition}
We shall call $\Gr=(U,V,E)$ the {\em  prime master-graph}.
In all bipartite graphs we consider in this text, the two sets
$U$ and $V$ of vertices remain constant, while the set 
of edges will vary through subsets of $E$ as defined in 
Definition \ref{d:bipartite}.%

\msk
\noindent
\subsection{Geometric Properties of the Master-Graph}
\label{sss:mastergraph}\index{Master-Graph}
The prime mas\-ter-graph can be drawn and helps in forming
a good intuition of the combinatorics involved.

\msk\index{graph!master}\index{master graph}
\index{vertex!upper}\index{vertex!lower}
\index{edge!sloping}

$\bullet$\quad  The set of vertices $U\cup V$ of the master-graph is
naturally contained in $\R^2=\R\times \R$, and so we can ``draw'' it
quite naturally. 

The elements in $\N\times \{1\}$ are called the {\em upper vertices},
those in $\N\times\{0\}$ the {\em lower vertices}

\msk

$\bullet$\quad  The edges $e=((n,1),(n,0))$, $n\in\N$ are called {\em vertical}.
All other edges $e=((m,1),(n,0))$ are called {\em sloping}. 
Because of $p_m|(q_n-1)$
they are sloping
``from left-above to right-below''. There are only vertical and sloping edges.
We call vertices $u=(m,1)$ and $v=(n,0)$ {\em connected} iff $e=(v,u)\in E$,
i.e., if $v$ is the upper vertex (end-point) of the edge $e$ and $v$ is the
lower vertex (end-point) of $e$.

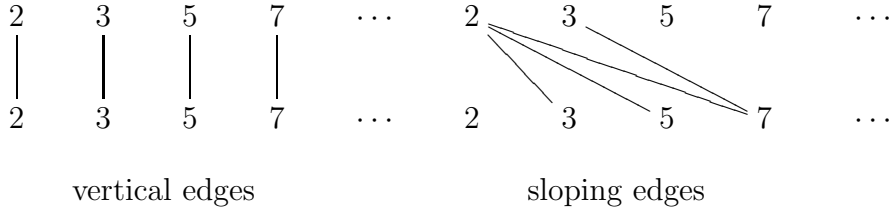
$$
\begin{tabular}{lcl}
\xymatrix@C=.7cm{2\ar@{-}[d]&3\ar@{-}[d]&5\ar@{-}[d]&7\ar@{-}[d]&\cdots\\
               2&3&5&7&\cdots}& & 
\xymatrix{2\ar@{-}[dr]\ar@{-}[drr]\ar@{-}[drrr]&3\ar@{-}[drr]&5&7&\cdots\\
               2&3&5&7&\cdots}\\
               &&\\
 \text{\ \ \ \ \ \ \ vertical edges} & &\text{\ \ \ \ \ \ \ sloping edges} 
               \end{tabular}
               $$
\msk

$\bullet$\quad  Each lower vertex $(n,0)$ is the endpoint of one vertical and
{\em finitely many} sloping edges.  It is connected to an upper vertex $(m,1)$
iff $p_m|(q_n-1)$.\\[1.5mm]
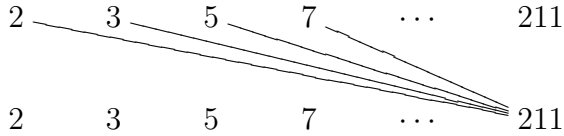
\begin{tabular}{l}
\xymatrix{2\ar@{-}[rrrrrd]&3\ar@{-}[rrrrd]&5\ar@{-}[rrrd]&7\ar@{-}[rrd]&\cdots&211\\
                2&3&5&7&\cdots&211&&&&} \\
                \\
{\small Edges in $U$ connected to the lower edge with label $211$ in $V$.  } 
\end{tabular}

\msk

$\bullet$\quad  Each upper vertex $(m,1)$ is connected to infinitely many
lower vertices $(n,0)$, namely, 
all those for which $p_m|(q_n-1)$, that is, 
for which there is a natural number $k$ for which $q_n=kp_m+1$. Indeed,
Dirichlet's Prime Number Theorem says: {\em Every arithmetic progression 
of the form $\{ka+b:k\in\N\}$ with $a$ and $b$ relatively prime,
contains infinitely many primes.} \index{Prime Number Theorem}
\index{arithmetic progression}

\begin{definition} \label{d: master} \rm  Let $p$ and $q$ be any primes, 
say, $p=p_m$ and $q=q_n$.
Then \index{master-graph!funnel}\index{master-graph!cone}
\index{funnel}\index{cone}
$$\cE_p=\{e: e=((m,1),(n,0))\in E\mbox{ such that }p|(q_n-1)\},$$
 the set of all edges emanating downwards from the vertex $(m,1)\in U$
will be called the {\em cone peaking at $p$}.
This cone contains infinitely many edges while
$$\cF_q=\{e: e=((m,1),(n,0))\in E\mbox{ such that }p_m|(q-1)\},$$
the set of edges ending below in the vertex $(n,0)\in V$,
called the {\em funnel pointing to $q$}, contains only
finitely many edges.
\end{definition}

\bsk

\begin{figure}
\begin{tikzpicture}  
\coordinate(u2) at (0,0);
\coordinate(u3) at (2,0);
\coordinate(u5) at (4,0);
\coordinate(u7) at (5.5,0); 
\coordinate(u11) at (7.5,0); 
\coordinate(u13) at (9,0); 

\coordinate(o2) at (0,4);
\coordinate(o3) at (2,4);
\coordinate(o5) at (4,4);
\coordinate(o7) at (5.5,4); 
\coordinate(o11) at (7.5,4); 
\coordinate(o13) at (9,4); 

\fill (u2) circle (2pt); \node at (u2)[below] {2};  
\fill (u3) circle (2pt); \node at (u3)[below] {3};
\fill (u5) circle (2pt); \node at (u5)[below] {5}; 
\fill (u7) circle (2pt); \node at (u7)[below] {7};  
\fill (u11) circle (2pt); \node at (u11)[below] {11};
\fill (u13) circle (2pt); \node at (u13)[below] {13};

\fill (o2) circle (2pt); \node at (o2)[above] {2};  
\fill (o3) circle (2pt); \node at (o3)[above] {3};
\fill (o5) circle (2pt); \node at (o5)[above] {5}; 
\fill (o7) circle (2pt); \node at (o7)[above] {7};  
\fill (o11) circle (2pt); \node at (o11)[above] {11};
\fill (o13) circle (2pt); \node at (o13)[above] {13}; 

\draw[thin]  (o2) --(u3);
\draw[thin]  (o2) --(u5);
\draw[thin]  (o2) --(u7);
\draw[thin]  (o2) --(u11);
\draw[thin]  (o2) --(u13);
\draw[thick]  (o3) --(u7);
\draw[thick] (o5) --(u11);
\draw[thick] (o3) --(u13);

\draw[thin] (o2)--(u2);
\draw[thin] (o3)--(u3);
\draw[thin] (o5)--(u5);
\draw[thin] (o7)--(u7);
\draw[thin] (o11)--(u11);
\draw[thin] (o13)--(u13);
\end{tikzpicture}
\caption{The initial part of the master graph.}
\end{figure}
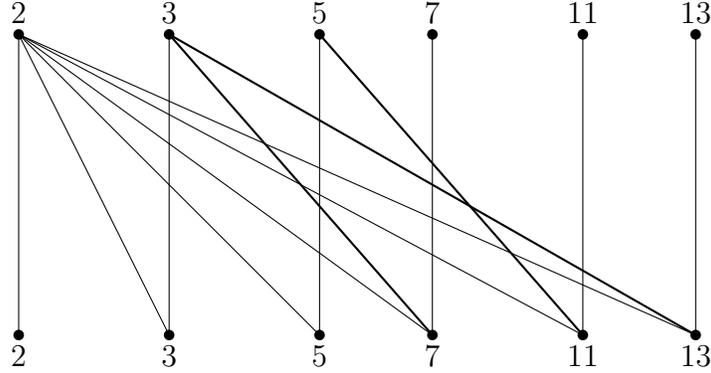
\ssk

All applications of prime graphs which we use in the structure theory
of \na\ groups 
are subgraphs of this master graph.
\bsk

\rm
Since for any periodic \lca\ group $A$ 
we have a canonical surjective morphism 
$\mu\colon {\til\Z}^\times\to \SAut(A)$
we need explicit information on the primary structure 
-- or $p$-Sylow structure -- of
${\til\Z}^\times$.   We are now going to describe this structure 
in additive notation in terms of the prime master-graph
$\Gr=(U,V,E )$.

\msk

Let $e=((m,1),(n,0))\in E$ be an edge in the master-graph.

Case 1. $m=n$. Then we set $\SS_e=\Z_{p_m}$.

\msk

Case 2.  $m<n$. Then $p_m|(q_n-1)$. Assume that the $q_n-1=p_m^{k(e)}s(e)$
with $s(e)$ relatively prime to $p_e\defi p_m$. 
Then we set $\SS_e=\Z(p_m^{k(e)})$

For the following proposition we recall 
$${\til\Z}^\times=\prod_{q\in\pi}\Z_q^\times.$$
Since $\Z_q^\times$ is not a \psyl q subgroup, this is {\em not}
the $q$-primary decomposition of $\til\Z^\times$. That decomposition
we describe now:

\begin{proposition} \quad {\rm(The  Sylow Structure of $\til\Z^\times$)}\quad 
\label{p:Z_q-times} Let $p, q\in\pi$ be  primes. 
Then \index{Sylow structure!of $\til\Z^\times$}
\begin{enumerate}[\rm(i)]
\item The structure of $\Z_q^\times$ (in additive notation) is
$$\prod_{e\in\cF_q}\SS_e
   =\Z_q\times\prod_{e\in\cF_q,\ {\rm sloping}}\Z(p_e^{k(e)}).$$
\item The $p$-primary component or
\psyl p subgroup of ${\til\Z}^\times$

$$({\til\Z}^\times)_p=\prod_{e\in\cE_p}(\Z_{q_e}^\times)_{p_e}$$
 is (in additive notation) 
$$({\til\Z}^\times)_p=\prod_{e\in\cE_p}\SS_e
  =\Z_p\times\prod_{e\in\cE_p,\ {\rm sloping}}\Z(p^{k(e)}).$$
\end{enumerate} 
\end{proposition}

\bsk

\subsection{The Structure of the Invertible Scalar Multiplications of an Abelian Group
 and their  Prime Graph}
\label{sss:A-master}\rm
Now let $A$ be a periodic \lca\ group; the Sylow structure of
$\SAut(A)$ is now easily discussed: The quotient morphism
$\mu\colon {\til \Z}^\times\to \SAut(A)$, preserving the Sylow structures, 
and the structure of $\SAut(A)$ described
so far in Theorem \ref{th:muk3-th2} allow a precise description
of the Sylow structure of $\SAut(A)$.

We associate with $A$ the bipartite graph $\Gr(A)=(U,V,E(A))$
with $U$ and $V$ as in the master-graph and with 
$$E(A)=\{e\in E: e=((m,1),(n,0))\mbox{ such that }\SAut(A_{q_n})
\ne\{\id_A\}\},$$
and we define

$\cE_p=\{e\in E(A): e=((m,1),(n,0))\in e(A)\mbox{ such that }p|(q_n-1)\}$,
the set of all edges in $\Gr(A)$ ending at the vertex $(n,0)\in V$
such that $\SAut(A_{q_n})$ is nontrivial, and

$\cF_q=\{e\in E(A): e=((m,1),(n,0))\in E(A)\mbox{ such that }p_m|(q-1)\}$,
the set of all edges in $\Gr(A)$ ending at the vertex $(n,0)\in V$ with
$q_n=q$ such that $\SAut(A_q)$ is nontrivial.

We recall  that for  each $q$-primary
component $A_q$, the ring of scalars $\SAut(A_q)$ is either
cyclic of order $q^r$, the exponent of $A_q$, if it is finite,
and is $\cong\Z_q$ otherwise. Thus its $q$-primary component is
$$\cong\begin{cases}\Z(q^{r-1})&%
       \mbox{if the exponent of $A_q$ is finite}\\ 
  \Z_q&\mbox{otherwise}.\end{cases}$$ 

Accordingly we define, for each
edge $e=((m,1),(n,0))\in E(A)$ in the graph $\Gr(A)$

$$\SS_e(A)=\begin{cases} \Z(q_m^{r-1})
&\mbox{if $m=n$, and $A_q$ has finite exponent $q^r$}\\
                \Z_{q_m}&\mbox{if $m=n$, and $A_q$ has infinite exponent,}\\
           \Z(p_m^{k(e)})&\mbox{if $m<n$}.\end{cases}$$
Then we have, analogously to Proposition \ref{p:Z_q-times}, the following
theorem, complementing 
Proposition \ref{p:Z_q-times} and Theorem \ref{th:muk3-th2}:

\begin{theorem}\quad{\rm(The Sylow Structure of $\SAut(A)$)}
    \label{c:Saut}
 Let $A$ be a periodic \lca\ group and
$\SAut(A)=\prod_{p\in\pi}\SAut(A)_p$ the $p$-primary decomposition
of the profinite group $\SAut(A)$.
 Then
\begin{enumerate}[\rm(i)]
\item The $p$-primary decomposition  of
$\SAut(A_q)$ is (additive notation assumed) 
$$\prod_{e\in\cF_q}\SAut(A_{q_e})_{p_e}\cong\prod_{e\in\cF_q}\SS_e(A)
  =\Z_q\times\prod_{e\in\cF_q,\ {\rm sloping}}\Z(p_e^{k(e)}).$$
\item The structure of the $p$-primary component $\SAut(A)_p$ 
of $\SAut(A)$ (in additive notation) is
$$\prod_{e\in\cE_p}(\Saut(A_{q_e})_{p_e}=\prod_{e\in\cE_p}\SS_e(A)
   =\Z_p\times\prod_{e\in\cE_p,\ {\rm sloping}}\Z(p_e^{k(e)}).$$
\end{enumerate} 
\index{Sylow structure!of $\SAut$}
\index{S Aut@$\SAut$!Sylow structure}
\end{theorem}

This theorem illustrates  the usefulness of the prime graph
$\Gr(A)$ which elucidates the fine structure of $\SAut(A)$.
In many instances, the prime graph is equally helpful in
the discussion of the Sylow structure of any periodic 
\na\ \lc\ group $G$ (see \cite{BOOK}, Section 10).

\section{The Prime Graph of a Near Abelian Group}

Staying with a periodic \na\ group $G$ which is $A$-nontrivial for a
base group $A$,  we investigate the interaction of the different 
Sylow subgroups in terms of the {\em prime graph} defined as a 
subgraph of the master-graph $\Gr=(U,V,E)$ of Definition \ref{d:bipartite}
as follows:

\begin{definition} \label{d:primegraphG} Let $G$ be a periodic \na\ $A$-nontrivial
\lc\ group with a base group $A$ and write $G=AZ(G)H$.

\centerline{%
A subgraph $\Gr_G=(U_G,V_G,E_G)$ of the master-graph $\Gr$}

\noindent  is called
the {\em prime graph of} $G$ provided the following conditions are satisfied:
\begin{enumerate}[\rm(i)]
\item 
We call $(m,1)$ the {\em upper $p$-vertex} iff $p=p_m$ and
$(n,0)$ the {\em lower $q$-vertex} iff $q=q_n$.
\item An edge $e=((m,1),(n,0))$ of the mastergraph is an edge in $E_G$
      if and only if $[H_{p_m},A_{q_m}]\ne\{1\}$. With $p=p_m$ and $q=q_n$ this 
      edge is written $e_{pq}$ and called an {\em edge leading from} $p$ to $q$.
 
\item $(m,1)$ is an upper vertex in $U_G$ iff $(G/C_G(A))_p\ne\{1\}$.

\item $(n,0)$ is a lower vertex in $V_G$ iff $A_p\ne\{1\}$.
\end{enumerate}

\end{definition}
We have a  much sharper conclusion:

\msk
\begin{theorem}\label{th:Apq} {\rm(Structure Theorem IV on Periodic
Near Abel\-ian Groups)} 
 Let $G$ be a periodic $A$-nontrivial \na\ group. 
Let $e_{pq}$ be an edge in $\Gr_G$.

Then we have the following conclusions (see \cite{BOOK}, {\rm Theorem 10.13}):
\begin{enumerate}[\rm(1)]
\item If $p\ne q$, that is, $e_{pq}$ is sloping, then $p\ne2$ and
$p|(q -1)$, but above all
\begin{enumerate}
\item[$(C_1)$] for $x\in G_p\setminus C_G(A_q)$
               the function  $a\mapsto [x,a]: A_q\to A_q$ is an 
               automorphism of $A_q$. In particular, $[x,A_q]=A_q$.
\item[$(C_2)$] $A_q\cap Z(G)=\{1\}$.
\end{enumerate}

\item If $p=q$, that is, $e_{pq}$ is vertical, then there is a unit 
      $s\in\Z_q^\times$ and a natural
     number $m\in\N$ such that $[x,a]=a^{q^ms}$ for all $a\in A_q$,
     that is, $[x,A_q]=A_q^{q^m}$, and
     $A_p\cap Z(G)$ has an exponent dividing $q^m$.
\end{enumerate}
\end{theorem}

The theorem gives an impression of the circumstances in which the intersection
$A_q\cap Z(G)_q$ can be nontrivial: The lower $q$-vertex has to be isolated
in the prime graph in such a case. 

\section{Application 1: The Classification of \\
Topologically  Quasihamiltonian Groups}

The following definition is due to {\sc F. K\"ummich} \cite{kummich1}:

\begin{definition} \rm A topological group $G$ is called 
{\em topologically quasihamiltonian} if $\overline{XY}=
\overline{YX}$ holds for any pair of closed subgroups $X$
and $Y$ of $G$. 
\end{definition}

This is equivalent to saying that 
$\overline{XY}$ is a closed subgroup whenever $X$ and $Y$ are
subgroups of $G$.

With the framework provided by near abelian locally compact groups,
it is possible to classify completely the class of \tqh\ \lc\ groups.
The classification proceeds in two steps: In a first step we classify
all \lc\ \tqh\ groups, and in a second step we classify all locally compact
\tqh\ groups in one fell swoop. 

For step 1 we need a definition:
 \begin{definition} \label{d:quater} \rm The groups $M_n$
\index{quaternion group!generalised}\index{generalised quaternion group}
defined by generators and relations for $n=2,3,\dots$
according to
$$M_n\defi\gp{a,b\mid b^{2^n}=1,\ b^{2^{n-1}}=a^2,\ bab^{-1}=a^{-1}}$$
are called {\em generalised quaternion groups}.
\end{definition}
These groups  also satisfy the relations
$$ a^4=1 \mbox{ and } [a,b]=a^2$$
and are fully characterised by the following explicit
construction: 
$$M_n\cong \frac{\Z(4)\rtimes\Z(2^n)}{\Delta},$$
where $\Z(2^n)$ acts on $\Z(4)$ by scalar multiplication with $\pm 1$
and where $\Delta$ is generated by $(s,t)$, $s=2+4\Z$ and 
$t=2^{n-1}+2^n\Z$. (Cf. \cite{hofrus}, Definition 5.8.)
We note that $M_2$ is (isomorphic to) 
the usual group of quaternions  
$Q_8=\{\pm1,\pm i,\pm j,\pm k\}$ of eight elements. 

Here is step 1:

\begin{theorem}\label{p:tqh-p}  
\index{topologically quasihamiltonian!$p$-group}
A \lc\  $p$-group $G$ is \tqh\ if and only if $G$ 
is \na\ with a base group $A$ and an inductively monothetic
$p$-group $G/A$ and 
at least one of the following statements holds:

\begin{enumerate}[{\rm (a)}]
\item $G$ is abelian.

\item There is a $p$-procyclic 
scaling group $H=\gen b$ such that $G=AH$
 and there is a natural number
$s\ge1$, respectively, $s\ge2$, if $p=2$, such that $a^b=a^{1+p^s}$
 for all $a$ in $A$. The group $G$ is
$A$-nontrivial.

\item $p=2$ and $G\cong A_2 \times M_n$,  where $A_2$ is an exponent $2$ \lc\ 
abelian group according, and
$M_n$  is the generalised quaternion group of order $2^{n+1}$. 
In this case, 
 $$A=A_2\times\<a\>\cong \Z(2)^{(I_1)}\times \Z(2)^{I_2}\times \Z(4)$$
with $a$ as in Definition \ref{d:quater} for suitable 
sets $I_1$ and $I_2$. The group $G$ is $A$-trivial.
\end{enumerate}
\end{theorem}

Next step 2 (see \cite{BOOK}, Theorem 13.9):
 
\begin{theorem}\label{th:periodic-tqh-gps}
 Let $G$ be a \lc\ periodic \tqh\  group. Then, for each 
$p\in\pi(G)$, the set of $p$-elements $G_p$ is a \tqh\ 
$p$-group, and there is a compact open subgroup $U_p$ in $G_p$
such that $G=G_{\nu(G)}$ is (up to isomorphism)
the local product of \tqh\ $p$-groups
$$G\cong\prod_{p\in\pi(G)}^{\rm loc}(G_p,U_p).$$
Conversely, every group isomorphic to such a local product is a \tqh\ group.
\end{theorem}

This theorem is proved with the aid of our Theorem \ref{th:nu-theorem}
For nonperiodic abelian groups we give an algorithmic description of \tqh\ \lc\
groups in \cite{BOOK}, Theorem 13.14.
Except for $p=2$ it turns out that topologically quasihamiltonian groups 
in the general locally compact domain are
the same thing as near abelian groups. For the exceptional compact
2-groups that are near abelian but fail to be \tqh\ see \cite{hofrus}. 
\msk

Theorem \ref{th:periodic-tqh-gps} can be visualised in terms of its
prime graph, that all connected components are either vertical edges and
its end points or are isolated vertices. If we allow ourselves the identification of
the connected components of the prime graph with the subgroups they represent,
we could reformulate 
Theorem \ref{th:periodic-tqh-gps} as follows:

\begin{theorem}\label{th:periodic-tqh-gps2}  Let $G$ be a \lc\ periodic \tqh\  group.
Then each connected component of the prime graph of $G$ represents a normal $p$-Sylow
subgroup and $G$ is a local direct product of these subgroups.
\end{theorem} 

\section{Application 2: The Classification of \\Topologically 
Modular  Groups} 

\newcommand{\tMg}{to\-po\-lo\-gi\-cal\-ly mo\-du\-lar group}

Recall that the closed subgroups of a topological group form a lattice
w.r.t.\ inclusion ``${\subseteq}$'' as partial order.

\begin{definition} \rm A topological group $G$ is called 
{\em topologically modular} if 
the lattice of closed subgroups is modular, that is, satisfies
the law $X\vee(Y\cap Z)= (X\vee Y)\cap Z$ whenever $X$ is a closed
subgroup of $Z$.
\end{definition}

This is equivalent to saying that the lattice of closed subgroups does not
contain a sublattice isomorphic to 

$$
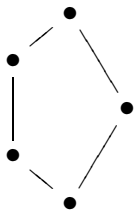
\xymatrix@R=2mm@C=3mm{& \bullet\ar@{-}[ld]\ar@{-}[rdd]& \\                                 
            \bullet\ar@{-}[dd]& & \\                                                         
                             & & \bullet\ar@{-}[ldd]\\                                       
            \bullet\ar@{-}[rd]& & \\                                                         
                 &\bullet&}$$

(See \cite{schmidt}, Theorem 2.1.2.)

It is instructive to spend some time on an example due to {\sc Mukhin} which shows 
that topologically modular groups can be tricky.

\begin{example}\label{ex:M}
Let $p$ be any prime and $I$ any infinite set (e.g. $I=\N$),
set $E\defi\Z(p)$, and define
define $G_j=E^2$, $C_j=\{0\}\times E$
 for all $j\in I$,  and set
$$G\defi E^{(I)}\times E^I\cong \prod_{j\in I}^{\rm loc}(G_j,C_j),$$
where we took the discrete topology on the direct sum $E^{(I)}$ and
the product topology on $E^I$.
 We shall identify $G$ with
$\prod_{j\in I}^{\rm loc}(G_i,C_i)$ and $E^I\times E^I$ with $(E^2)^I$.
The natural injection $\iota\colon G\to (E^2)^I=E^I\times E^I$
is continuous but is not an embedding, since it is not open
onto its image.

Let $D\defi\{(x,x):  x\in E\}\subseteq E^2$, and
$$\Delta =D^I=\{(x_j,x_j)_{j\in I}: x_j\in E\}\subseteq
(E^2)^I\cong E^I\times E^I$$ denote the respective diagonals.
Then $\Delta$ is a closed subgroup of $(E^2)^I$
and so $\iota^{-1}(\Delta)=D^{(I)}$ is a closed subgroup of $G$.
We shall denote it by $Y$. This is a noteworthy and perhaps
slightly unexpected fact in view of the density of $E^{(I)}$
in $E^I$. We verify as an exercise
 that the subgroup $Y$ is not only closed, but even discrete, since
$\iota(Y)$ meets  trivially every open subgroup $\{0\}\times E^K$ for
a  cofinite subset $K\subseteq I$.

Now the product $E^I$ is the 
projective limit of its finite partial products $E^F$ as $F$ ranges
through the directed set $\mathcal F$
of finite subsets $F$ of $I$. Accordingly,
$$G\cong E^{(I)}\times\lim_{F\in{\mathcal F}}E^F                                             
   \cong \lim_{F\in{\mathcal F}}(E^{(I)}\times E^F).$$

Let
$D_2=\{(x_j)_{j\in I}\in E^I: (\exists c\in E)(\forall j\in I)\, x_j=c\}$.
Now we consider the following subgroups of $G$:
$$\begin{matrix} X&\defi&E^{(I)}&{\times}&\{0\},\cr                                          
                 Z&\defi&E^{(I)}&{\times}&D_2\cr \end{matrix}$$
whence $X\subseteq Z$.
Then $X\vee Y=\overline{E^{(I)}\times E^{I}}=G$ and so
$(X\vee Y)\cap Z=Z$ on the one side, while  $Y\wedge Z= Y$ and so
$X\vee(Y\wedge Z)=X\vee Y=G$. Hence
$X\vee(Y\wedge Z)\ne (X\vee Y)\wedge Z$.
Therefore $G$ is a locally compact abelian nonmodular group.
\end{example}

The example shows that the  limit of a 
projective system
 of locally compact \tMg\ with  proper bonding maps need not be
a \tMg\ and that
a local product of a collection of finite abelian modular groups
may fail likewise to be a \tMg.

Locally compact abelian \tMg s were classified by {\sc Muk\-hin}
in \cite{muk2}. We now discuss the nonabelian situation.

A first step in the classification is the case of $p$-groups:

\begin{proposition}\label{p:p-tqh=p-tM} 
Let $G$ be a compactly ruled  $p$-group. Then the following statements
are equivalent:
\begin{itemize}
\item[(1)] $G$ is a topologically  modular group.

\item[(2)] $G$ is a \tqh\ group with  a base group $A$
that is a topologically modular \lca\ group.
\end{itemize}
\end{proposition}

In contrast, however, with topologically quasihamiltonian groups, the normal
$p$-Sylow subgroups are not the only building blocks of topologically modular 
groups.
There is one additional category of building blocks which in the discrete
situation were known since the pioneering work of Iwasawa in the forties of
the last century, see e.g. \cite{Iwa3}.

\subsection{Iwasawa (p,q)-Factors}

\begin{example} \label{ex:two-primes} For a prime $q$ let
{\bf A} be an additively written
\lca\ group of exponent $q$ and be either compact or discrete. 
Thus, algebraically, {\bf A} \
is a vector space over the field ${\rm GF}(q)$.

\ssk

Now let $p$ be a prime such that $p|(q-1)$. Then the multiplicative
group of GF$(q)$ contains a cyclic subgroup $Z$ of order $p$.
Let $C=\gen t$ be any $p$-procyclic group (that is,
$C\cong\Z(p^k)$ for some $k\in\N$ or $C\cong\Z_p$), and let
$\psi\colon C\to Z$ be an epimorphism. Then $C$
acts on {\bf A} via $r{*}a=\psi(r)\.a$. Since $Z$ is of order $p$,
the kernel of $\psi$ is an open subgroup of $C$ of index $p$.

Set $G={\bf A}\rtimes_\psi C$, the semidirect product for the action of
$C$ on {\bf A}. Then $A\defi {\bf A}\times\{1\}$ is a base subgroup of the
\na\ \lc\ group $G$, and  $H=\gen{(0,t)}=\{0\}\times C$ is a
procyclic scaling $p$-subgroup.

There are many maximal
$p$-subgroups of $G$, namely, each $\gen{(a,t)}$ for any
$a\in A$, and there is one unique maximal $q$-subgroup which is
normal, namely, $A$.

The simplest case arises when we take for $C$ the unique cyclic group
$S_p(Z)$  of $Z$ of order $p$, in which case we have
$G\cong {\bf A}\rtimes\Z(p)$ and the set of elements of order $p$ is
${\bf A}\times (\Z(p)\setminus\{0\})$ and the set of $q$-elements
is ${\bf A}\times\{0\}$.
\end{example}

The class of \lc\ \na\ \tqh\  groups described in
Example \ref{ex:two-primes} is relevant enough in our
classification to deserve a name:

\begin{definition} \label{d:iwasawa-factor}\rm
A \lc\ group $G$ which is isomorphic to a semidirect product
${\bf A}\rtimes_{\psi}C$ as described in Example \ref{ex:two-primes}
will be called an {\em Iwasawa $(p,q)$-factor}. The primes $p$ and $q$
are called the {\em primes of the factor $G$}.
\end{definition}

The prime graph $\mathcal G$ of an Iwasawa $(p,q)$-factor is one sloping
edge $e_{pq}$  with its endpoints.

We would like to see an abstract characterisation of a $(p,q)$-factor.
For the purpose of presenting one let us formulate some terminology
for an automorphic action $(h,a)\mapsto h\.a:H\times A \to A$
inducing a morphism $\alpha\colon H\to\Aut(A)$, $\alpha(h)(a)=h\.a$.
If $H/\ker\alpha$ is an abelian group of order $p$ for a prime
number $q$, we shall say that 

{\em the action of $H$ on $A$ is of order $p$.}

If $H$ is a subgroup of a group $G$ and $A$ is a normal subgroup of $G$,
then $H$ acts on $A$ via $h\.a=hah^{-1}$. If this action is of order $p$,
we say that

{\em $H$ induces an action of order $p$ on $A$.}

\begin{proposition}\label{p:iwasawa-pq} Let $p$ and $q$ be primes satisfying
$p|(q-1)$. 
A \na\ group $G$ is an Iwasawa $(p,q)$-factor if and only if it satisfies the
following conditions:
\begin{enumerate}[{\rm (a)}]
\item $A=G'$ is an abelian group of exponent $q$; it is either compact or discrete subgroup of $G$;
\item There is a scaling group $H$ which is a procyclic $p$-group; it induces an action of order $p$ on $A$.
\end{enumerate}
If these conditions are satisfied, then 
$G=A\rtimes H$ is a semidirect product
and $Z(G)=\{h^p: h\in H\}$.
\end{proposition}

The significance of the $(p,q)$-factors for our classification is
due to the following fact which requires a technical proof that is not
exactly short:

\begin{proposition}\label{p:Iwasawa-pq}
Let $G=AH$ be an Iwasawa $(p,q)$-factor and $A$ a topologically modular
abelian group. Then $G$ is a 
topologically modular group.
\end{proposition}

Since a nondegenerate $(p,q)$-factor does not meet the
criteria of a topologically quasihamiltonian \lc\ group 
in Theorem \ref{th:periodic-tqh-gps}, this allows us to remark a significant difference between
topologically quasihamiltonian and topologically modular groups: 

\begin{corollary} \label{c:qht-M}
 Any nondegenerate Iwasawa $(p,q)$-factor
provides a topologically modular group which is not 
topologically quasihamiltonian.
\end{corollary}

After a thorough discussion of compactly ruled topologically modular
groups, using much of the information accumulated on near abelian
groups we arrive at the following classification of periodic
locally compact topologically modular groups:

\begin{theorem}[The Main Theorem on Topologically Modular Groups]  
\label{th:tMg-summing-up} Let $G$ be a compactly ruled topologically
modular group. Then $\pi$ is a disjoint union of a set $J$ of 
sets $\sigma$ of prime numbers which are either empty, or singleton sets
$\sigma=\{p\}$ such that $G_\sigma$ is a normal $p$-Sylow subgroup and 
an Iwasawa $p$-factor,
or two element sets $\{p,q\}$ such that for $p<q$ the set $G_\sigma$ is
a normal $\sigma$-Sylow subgroup and an Iwasawa $(p,q)$-Factor,
such that 
$$ G=\prod_{\sigma\in J}^{\rm loc}(G_\sigma,C_\sigma)$$
for a family 
of compact open subgroups 
$C_\sigma\subseteq G_\sigma$. In particular, $G$ is a periodic \na\ \lc\
group.

Conversely, every near abelian locally compact $G$ 
of this form is a topologically modular locally
compact group.
\end{theorem}

We notice that, in the prime graph of $G$,
the Sylow $p$-subgroups $G_p$ constitute the connected components of either
isolated vertices or vertical edges with its endpoints, while the
Sylow subgroups $G_{\{p,q\}}$ which are Iwasawa $(p,q)$-components are
connected components consisting of sloping edges with their endpoints.
Moreover, every prime graph having such connected components can be realised
as the prime graph of a periodic locally compact topologically modular
group.

\section*{Acknowledgements}
We gratefully acknowledge Sidney A. Morris for his reading our text
and helping us to make the introduction to a complicated network of
mathematics more lucid. The second author looks back with
deep appreciation on at least
45 years of mathematical contacts about locally compact groups
with {\sc Herbert Heyer},
to whom this survey is dedicated.

\end{document}